\documentclass{amsart}
\usepackage{graphicx}
\usepackage{latexsym}
\usepackage{hyperref}

\usepackage{amsfonts}
\usepackage[all]{xy}

\numberwithin{equation}{section}
\newtheorem{thm}{Theorem}[section]

\theoremstyle{definition}

\numberwithin{equation}{section}
\theoremstyle{remark}

\newcommand{\be}{\begin{eqnarray}}
\newcommand{\en}{\end{eqnarray}}
\newcommand{\no}{\nonumber}
\newcommand{\laa}{\lambda}
\newcommand{\va}{\varsigma}

\newcommand{\ov}{\overline}

\newcommand{\ma}{\mathcal{E}}

\newcommand{\fr}{\frac}
\newcommand{\pa}{\partial}

\newcommand{\D}{\Delta}

\newcommand{\na}{\nabla}
\newcommand{\lan}{\langle}
\newcommand{\ran}{\rangle}
\newcommand{\lam}{\lambda}

\newcommand{\vs}{\vskip0.3cm}

\begin{document}
\title[Sharp  Bounds for the First  Eigenvalues of  the Bi-Laplace Operator] {Sharp Lower Bounds for the First  Eigenvalues of  the Bi-Laplace Operator}
\author[Qiaoling Wang] {Qiaoling Wang$^*$}
\
\thanks{$^{*}$Partially supported by CNPq-BR grant number: 307089/2014-2}
\author[Changyu Xia] {Changyu Xia$^{\dag}$}
\thanks{$^{\dag}$Partially supported by CNPq-BR grant number: 306146/2014-2.}
\thanks{$^{\dag}$Corresponding author}

\maketitle

\begin{abstract}

 We obtain sharp lower bounds for the
first eigenvalue of four types of eigenvalue problem defined by the
bi-Laplace operator on compact manifolds with boundary and determine
all the eigenvalues and the corresponding eigenfunctions of a
Wentzell-type bi-Laplace problem on Euclidean balls.
\end{abstract}
\footnotetext{{\it Key words and phrases}:
 \hspace*{2ex}Key words and phrases: Eigenvalue,
Bi-Laplace Operator,  Ricci curvature, Steklov Problem,
 Euclidean ball, hemisphere. }

\renewcommand{\sectionmark}[1]{}

\section{Introduction and the main results}

The classical Lichnerowicz-Obata theorem is stated as follows (Cf.
\cite{Chav84,Lich,obata}): \vskip0.2cm {\bf Theorem A}
(Lichnerowicz-Obata). {\it Let $M$ be an $n$-dimensional complete
Riemannian manifold with Ricci curvature bounded below by
$(n-1)k>0$, then the first non-zero eigenvalue of the Laplacian of
$M$ satisfies \be\label{thbeq1} \eta \geq n\kappa,\en and equality
holds if and only if $M$ is isometric Euclidean $n$-sphere of radius
$1/\sqrt{\kappa}.$ } \vskip0.2cm
 When $M$ has nonnegative Ricci
curvature, Li-Yau \cite{LY} proved that \be \label{ly}
\eta_1\geq\fr{\pi^2}{2D^2}, \en where $D$ is the diameter of $M$.
Later, Zhong-Yang  improved this estimate to an optimal lower bound
\cite{ZY}: \be \label{zy1}\lam\geq \fr{\pi^2}{D^2}. \en Moreover,
Hang-Wang \cite{HW} showed  that if the equality in (\ref{zy1})
holds then M must be isometric to a circle of radius $D/\pi$.

On the other hand,  in 1977, Reilly obtained a sharp lower for the
first Dirichlet eigenvalue of the Laplacian of compact manifolds
with boundary and positive Ricci curvature \cite{RE}. \vskip0.2cm
 \vskip0.2cm{\bf Theorem B} (Reilly).
 {\it Let $M$ be an
$n(\geq 2)$-dimensional compact  Riemannian manifold with Ricci
curvature bounded below by $(n-1)\kappa>0$ and  boundary.  If the
mean curvature of $\pa M$ is nonnegative then the first Dirichlet
eigenvalue of the Laplacian of $M$ satisfies \be\label{re1}
\lam_1\geq n\kappa,\en  and equality holds if and only $M$ is
isometric to an $n$-dimensional Euclidean hemisphere of radius
$1/\sqrt{\kappa}$. } \vskip0.2cm
 For the Neumann boundary case a
similar result has been proven by Escobar \cite{Esc90} and Xia
\cite{X1} independently. They showed that if $M$ is an $n(\geq
2)$-dimensional compact Riemannian manifold with Ricci curvature
bounded below by $(n-1)\kappa>0$ and convex boundary, then the first
nonzero eigenvalue of the Laplace of $M$ with  Neumann boundary
condition satisfies \be\label{EX} \mu_1 \geq n\kappa,\en
 and equality holds if and only if $M$ is isometric
to and $n$-dimensional Euclidean hemisphere of curvature $\kappa$.

A natural question related to the above results is to study similar
rigidity phenomenon  for  other eigenvalue problems on Riemannian
manifolds. In the first part of this paper, we consider two kinds of
eigenvalue problem of the bi-Laplace operator on compact Riemannian
manifolds with boundary and obtain sharp lower bound for the first
eigenvalue of them. The point of interest is that no assumption on
the boundary is made.

\begin{thm}\label{th5}
 Let $M$ be an $n(\geq 2)$-dimensional compact Riemannian manifold with boundary.
 Assume that the Ricci
 curvature of $M$ is bounded below by $(n-1)k\geq 0$.
 Denote by $\lambda_1$  the first  Dirichlet eigenvalue  of the Laplacian of $M$
 and  $\Gamma_1$  the first eigenvalue of the following
 problem :
\be \label{th5eq1} \left\{\begin{array}{l}
  \Delta^2 u= \Gamma  u   \ \ {\rm in \ \ } M, \\
 u= \fr{\pa^2 u}{\pa \nu^2}=0 \ \ {\rm on \ \ } \pa M.
\end{array}\right.
\en Then \be\label{th5eq2} \Gamma_1\geq
\lambda_1\left(\fr{\lambda_1}n+(n-1)\kappa\right), \en with equality
holding if and only if $\kappa >0$ and $M $ is isometric to an
$n$-dimensional Euclidean hemisphere of radius $ 1/\sqrt{\kappa}$.
\end{thm}

\begin{thm}\label{th6}
Let $M$ be an $n(\geq 2)$-dimensional compact Riemannian manifold
with boundary.
 Assume that the Ricci
 curvature of $M$ is bounded below by $(n-1)k$.
 Denote by $\lambda_1$  the first  Dirichlet eigenvalue  of the Laplacian of $M$
and $\Lambda_1$  the first eigenvalue of the
 problem
\be \label{th6eq1} \left\{\begin{array}{l}
  \Delta^2 u= -\Lambda \Delta u   \ \ {\rm in \ \ } M, \\
 u= \fr{\pa^2 u}{\pa \nu^2}=0 \ \ {\rm on \ \ } \pa M.
\end{array}\right.
\en Then \be \label{th6eq2} \Lambda_1\geq
\fr{\lambda_1}n+(n-1)\kappa, \en
 with equality holding if and only if
$\kappa>0$ and $M$ is isometric to an $n$-dimensional Euclidean
hemisphere of radius $1/\sqrt{\kappa}$.
\end{thm}
The eigenvalue problems (\ref{th5eq1}) and (\ref{th6eq1}) should be
compared with the clamped plate problem and the buckling problem,
respectively. The later two ones are as follows: \be \label{clameq1}
\left\{\begin{array}{l}
  \Delta^2 u= \Gamma  u   \ \ {\rm in \ \ } M, \\
 u= \fr{\pa u}{\pa \nu}=0 \ \ {\rm on \ \ } \pa M.
\end{array}\right.
\en \be \label{buckeq1} \left\{\begin{array}{l}
  \Delta^2 u= -\Lambda \Delta u   \ \ {\rm in \ \ } M, \\
 u= \fr{\pa u}{\pa \nu}=0 \ \ {\rm on \ \ } \pa M.
\end{array}\right.
\en

Theorems (\ref{th5}) and (\ref{th6}) shows that the first eigenvalue
of (\ref{th5eq1}) (or (\ref{th6eq1})) is closely related to  the
first Dirichlet eigenvalue of the Laplacian. Many important results
have been obtained for the eigenvalues of the problems
(\ref{clameq1}) and (\ref{buckeq1}) (cf. \cite{ABe}, \cite{ABL},
\cite{ABu}, \cite{AL}, \cite{CY1}, \cite{CY2}, \cite{CY3},
\cite{CY4}, \cite{H},  \cite{JLWX1}, \cite{JLWX2}, \cite{N},
\cite{WX1}, \cite{WX2}, \cite{WX3}, etc). It would be interesting to
characterize the  hemisphere by the first non-zero eigenvalue of the
problem (\ref{clameq1}) or (\ref{buckeq1}).

In the second part of this paper,  we study Steklov eigenvalue
problems of the bi-Laplace operator. The eigenvalue problems we are
interested in are as follows: \be
\label{eq1.5}\left\{\begin{array}{l}
 \Delta^2 u= 0  \ \ {\rm in \ \ } M, \\
u= \D u  -p\fr{\pa u}{\pa \nu} =0 \ \ {\rm on \ \ } \pa M,
\end{array}\right.
\en \be \label{eq1.6}\left\{\begin{array}{l}
 \Delta^2 u= 0  \ \ {\rm in \ \ } M, \\
u= \fr{\pa^2 u}{\pa \nu^2}  -q\fr{\pa u}{\pa \nu} =0 \ \ {\rm on \ \
} \pa M,
\end{array}\right.
\en \be \label{eq1.7}\left\{\begin{array}{l}
 \Delta^2 u= 0  \ \ {\rm in \ \ } M, \\
\fr{\pa u}{\pa \nu} = \fr{\pa(\D u) }{\pa \nu}  + \xi u =0 \ \ {\rm
on \ \ } \pa M,
\end{array}\right.
\en
 and \be \label{eq1.8}\left\{\begin{array}{l}
 \Delta^2 u= 0  \ \ {\rm in \ \ } M, \\
\fr{\pa u}{\pa \nu} = \fr{\pa(\D u) }{\pa \nu} +\beta \overline{\D}
u + \varsigma u =0 \ \ {\rm on \ \ } \pa M,
\end{array}\right.
\en where $\beta$ is a nonnegative constant and $\overline{\D}$ is
the Laplace operator with respect to the induced metric of $\pa M$.

Elliptic problems with parameters in the boundary conditions are
called Steklov problems from their first appearance in \cite{S}.
Problem (\ref{eq1.5}) was considered by Kuttler \cite{K1} and Payne
\cite{P2} who studied the isoperimetric properties of the first
eigenvalue $p_1$ which is the sharp constant for $L^2$ a priori
estimates for solutions of the (second order) Laplace equation under
nonhomogeneous Dirichlet boundary conditions (cf. \cite{K2,KS1,KS2}
). The whole spectrum of  (\ref{eq1.5}) was studied in
\cite{FGW1,LIU1,LIU2} where one can also find a physical
interpretation of $p_1$. We refer to \cite{BGM,BFG,BP,GS,RS,WX4} for
some further developments about $p_1$. One can  see that $p_1$ is
positive and given by \be\label{eq1.9} p_1=\min_{w|_{\pa M}=0, w\neq
const.}\fr{\int_M (\D w)^2}{\int_{\pa M}\left(\fr{\pa w}{\pa
\nu}\right)^2}. \en

The problem (\ref{eq1.6}) is a natural Steklov problem and  is
equivalent to (\ref{eq1.5}) when the mean curvature of $\pa M$ is
constant. We have a sharp relation between the first eigenvalues of
(\ref{eq1.5}) and (\ref{eq1.6}).
\begin{thm}\label{th7} Let $M$ be an
$n$-dimensional compact Riemannian manifold with boundary and
non-negative Ricci curvature. Denote by $p_1$ and $q_1$ the first
eigenvalue of the problems (\ref{eq1.5}) and (\ref{eq1.6}),
respectively. Then we have \be \label{th7eq1} q_1\geq \fr{p_1}n,\en
with equality holding if and only if $M$ is  isometric to a ball  in
$ \mathbb{R}^n$.
\end{thm}
The next result is  a sharp lower bound for $p_1$.
\begin{thm}\label{th8}
 Let $M$ be a compact Riemannian manifold with  Ricci curvature bounded  below by
$-(n-1)\kappa$ for some constant $\kappa \geq 0$. Denote by
$\lambda_1$ the first Dirichlet eigenvalue of the Laplacian of $ M$
and let $p_1$  be the first eigenvalue of the problems
(\ref{eq1.5}). Suppose that  the mean curvature of $\pa M$ is
bounded below by a  positive constant $c$. Then we have \be
\label{th8eq1} p_1\geq \fr{nc\lambda_1}{n\kappa +\lambda_1} \en with
equality holding if and only if $\kappa=0$ and $M$ is isometric to a
Euclidean $n$-ball of radius $1/c$.
\end{thm}
The problem (\ref{eq1.7}) was first studied in (\cite{KS1}) where
some estimates for the first non-zero eigenvalue $\xi_1$ were
obtained. When $M$ is an Euclidean ball, all the eigenvalues of the
problem (\ref{eq1.6}) have been  recently obtained in \cite{XW}.
Also, the authors  proved an isoperimetric upper bound for $\xi_1$
when $M$ is a bounded domain in $\mathbb{R}^n$. The Rayleigh-Ritz
formula for $\xi_1$ is: \be\label{pb1.7} \xi_1=\underset{\int_{\pa
M}u=0, \pa_{\nu}u|_{\pa M}=0}{\underset{0\neq u\in
H^2(M)}{\min}}\frac{\int_{M} (\D u)^2} {\int_{\pa M} u^2}. \en The
problem (\ref{eq1.8}) is a so called Wentzell problem for the
bi-laplace operator which is motivated by (\ref{eq1.7}) and the
following Wentzell-Laplace problem:
\be\label{wen}\left\{\begin{array}{l} \Delta u
=0 \ \ \ \ \ \ \ \ \ \ \ \ \ \ \ \ \ \ {\rm in  \ }\ M,\\
-\beta\ov{\Delta} u+\pa_{\nu} u=  \laa  u\ \ \ {\rm on  \ } \pa M,
\end{array}\right.
\en where $\beta$ is a given non-negative  number. The  problem
(\ref{wen}) has been studied recently, in \cite{DK}, \cite{XW}, etc.

The first non-zero eigenvalue of (\ref{eq1.8}) can be characterized
as \be \varsigma_{1,\beta}=\underset{\int_{\pa M}u=0,
\pa_{\nu}u|_{\pa M}=0}{\underset{0\neq u\in H^2(M)}{\min}}\fr{\int_M
(\D w)^2 +\beta\int_{\pa M} |\overline{\na} w|^2}{\int_{\pa M} w^2}.
\en From (\ref{pb1.7}), one can see that if $\beta>0,\ \xi_1$ is the
first non-zero eigenvalue of the Steklov problem (\ref{eq1.7}) and
$\lambda_1$ the first non-zero eigenvalue of the Laplacian  of $\pa
M$, then we have \be\label{eqxi} \varsigma_{1,\beta}\geq \xi_1+\beta
\lambda_1, \en with equality holding if and only if any
eigenfunction $f$ corresponding to $\varsigma_{1, \beta}$ is an
eigenfunction corresponding to $\xi_1$ and $f|_{\pa M} $ is an
eigenfunction corresponding to $\lambda_1$.

For each $k=0, 1,\cdots,$ let $\mathcal{D}_k$ be the space of
harmonic homogeneous polynomials in $\mathbb{R}^n$ of degree $k$ and
denote by $\mu_k$ the dimension of $\mathcal{D}_k$. We refer to
\cite{ABR} for the basic properties of  $\mathcal{D}_k$ and $\mu_k$.
In particular, we have \be\no & &\mathcal{D}_0={\rm span}\{1\}, \ \
\ \mu_0=1,\\ \no &
&\mathcal{D}_1={\rm span}\{x_i, i=1,\cdots, n\}, \  \ \ \mu_1=n,\\
\no & &\mathcal{D}_2={\rm span}\{x_ix_j, x_1^2-x_h^2,  1\leq i<j\leq
n, h=2,\cdots, n\}, \ \ \ \ \mu_2=\fr{n^2+n-2}2. \en

In the next result we determine explicitly all the eigenvalues of
(\ref{eq1.8}) and the corresponding eigenfunctions  when
$M=\mathbb{B}$(the unit ball with center at the origin in $\mathbb
{R}^n$).

\begin{thm}\label{th9} If $n\geq 2$ and $M =\mathbb{B}$, then we have

 i) the eigenvalues of (\ref{eq1.8}) are $\va_k=k^2(n+2k)+\beta k(k+n-2),\ k=0, 1, 2, \cdots $;

 ii) for all $k=0, 1, 2,\cdots,$ the multiplicity of $\va_k$ is $\mu_k$;

 iii) for all $k=0, 1, 2,\cdots,$ the eigenspace corresponding to $\va_k$ is given by
 \be \mathbf{E}_k=\{ -2w_k+k(|x|^2-1)w_k|\ \  w_k\in {\mathcal D}_k\}.
\en
\end{thm}
Our last result is a lower bound for $\va_1$.
\begin{thm}\label{th10}
Let $M$ be an $n$-dimensional compact Riemannian manifold with
boundary and  Ricci curvature bounded below by $-(n-1)\kappa, \
\kappa\geq 0$. Assume that the principal curvatures of $\pa M$ are
bounded below by a positive constant $c$ and denote by  $\va_1$ the
first eigenvalue of the problems and (\ref{eq1.8}). Then we have \be
\label{th10eq1} \va_1> \fr{nc\lambda_1 \mu_1}{(n-1)(\mu_1+n\kappa
)}+\beta\lambda_1, \en where $\mu_1$ and $\lambda_1$ are the first
nonzero Neumann eigenvalue of the Laplacian of $M$ and the first
nonzero eigenvalue of the Laplacian of $\pa M$, respectively.
\end{thm}

\section{Proof of Theorems \ref{th5} and \ref{th6} }

In this section, we shall prove Theorems \ref{th5} and \ref{th6}.
 Before doing this, we first recall Reilly's formula which will be used later. Let $M$
be an $n$-dimensional compact manifold with boundary. We will often
write $\langle, \rangle$ the Riemannian metric on $M$ as well as
that induced on $\pa M$. Let $\na$ and $\D $ be the connection and
the Laplacian on $M$, respectively. Let $\nu$ be the unit outward
normal vector of $\pa M$. The shape operator of $\pa M$ is given by
$S(X)=\na_X \nu$ and the second fundamental form of $\pa M$ is
defined as $II(X, Y)=\langle S(X), Y\rangle$, here $X, Y\in T \pa
M$. The eigenvalues of $S$ are called the principal curvatures of
$\pa M$ and the mean curvature $H$ of $\pa M$ is given by $H=\fr
1{n-1} {\rm tr\ } S$, here ${\rm tr\ } S$ denotes the trace of $S$.
 For a smooth function
$f$ defined on $M$, the following identity holds \cite{RE} if
$h=\left.\fr{\pa }{\pa \nu}f\right|_{\pa M}$, $z=f|_{\pa M}$ and
${\rm Ric}$ denotes the Ricci tensor of $M$: \be & & \int_M
\left((\D f)^2-|\na^2 f|^2-{\rm Ric}(\na f, \na f)\right)\\ \no &=&
\int_{\pa M}\left( ((n-1)Hh+2\overline{\D}z)h + II(\overline{\na}z,
\overline{\na}z)\right). \en Here $\na^2 f$ is the Hessian of $f$;
$\ov{\D}$ and $\ov{\na}$  represent the Laplacian and the gradient
on $\pa M$ with respect to the induced metric on $\pa M$,
respectively.
 \vskip0.2cm

{\it Proof of Theorem \ref{th5}.} Let $f$ be an eigenfunction of the
problem (\ref{th5eq1}) corresponding to the first eigenvalue
$\Gamma_1$. That is,
 \begin{equation}\label{th5eq3}
\left \{ \aligned &\Delta^2 f= \Gamma_1  f   \ \ {\rm in \ \ } M, \\
 &f= \fr{\pa^2 f}{\pa \nu^2}=0 \ \ {\rm on \ \ } \pa M.
 \endaligned \right.
 \end{equation}
It  follows from the divergence theorem that \be \label{th5eq4}
\Gamma_1\int_M
f^2&=& \int_M f\D^2 f\\ \no &=&-\int_{M}\lan \na f, \na(\D f)\ran \\
\no &=& \int_M (\D f)^2 -\int_{\pa M} h\D f, \en where
$h=\left.\fr{\pa f}{\pa \nu}\right|_{\pa M}.$ Since $f|_{\pa
M}=\left. \fr{\pa^2 u}{\pa \nu^2}\right|_{\pa M}=0$, we have \be
\label{th5eq5} \D f|_{\pa M} = (n-1) H h. \en From Reilly's formula,
we infer \be \label{th5eq6}
 \int_M \left((\D f)^2-|\na^2 f|^2\right)= \int_M{\rm Ric}(\na f, \na f)+(n-1)\int_{\pa M}Hh^2. \en
Combining (\ref{th5eq4})-(\ref{th5eq6}) and using ${\rm Ric}(\na f,
\na f)\geq (n-1)\kappa |\na f|^2$ , we get \be \label{th5eq7}
\Gamma_1&=&\fr{\int_M (|\na^2 f|^2+{\rm Ric}(\na f, \na f)}{\int_M
f^2}\\ \no &\geq& \fr{\int_M (|\na^2 f|^2+ (n-1)\kappa|\na
f|^2)}{\int_M f^2}. \en The Schwarz inequality implies that \be
\label{th5eq8} |\na^2 f|^2\geq \fr 1n (\D f)^2, \en with equality
holding if and only if \be\label{th5eq9} \na^2 f =\fr{\D f}n \lan ,
\ran. \en
Therefore, we have \be  \Gamma_1\geq\fr{\int_M \left(\fr 1n (\D
f)^2+ (n-1)\kappa|\na f|^2\right)}{\int_M f^2} , \en with equality
holding if and only (\ref{th5eq9}) holds and \be\no {\rm Ric}(\na f,
\na f)= (n-1)\kappa |\na f|^2 \ \ {\rm on}\  M. \en

On the other hand, since $f$ is not a zero function which vanishes
on $\pa M$, we know that \be \label{th5eq10} \int_M (\D f)^2\geq
\lambda_1\int_M |\na f|^2\geq \lambda_1^2 \int_M f^2 \en with
equality holding if and only if $f$ is a first eigenfunction of the
Dirichlet Laplacian of $M$. Thus we conclude that \be
 \Gamma_1\geq
\lambda_1\left(\fr{\lambda_1}n+(n-1)\kappa\right).\en Assume now
that \be\label{th5eq11} \Gamma_1=
\lambda_1\left(\fr{\lambda_1}n+(n-1)\kappa\right).\en  In this case,
(\ref{th5eq10}) should take equality sign which implies that $f$ is
a first eigenfunction corresponding to the the first eigenvalue
$\lambda_1$ of the Dirichlet Laplacian of $M$. Consequently, we have
 \be
  \Delta f= -\lambda_1  f   \ \ {\rm in \ \ } M
\en
It then follows that
 \be
  \Delta^2 f= -\lambda_1\D  f =\lambda_1^2 f \ \ {\rm in \ \ M}
 \en
 which, combining with (\ref{th5eq3}) and (\ref{th5eq11}), gives $\lambda_1=n\kappa >0$. Since $f$ is an eigenfunction
 corresponding to the first Dirichlet eigenvalue of the Laplacian, it is well known that
 $f$ doest not change sign in $M$. Let us assume that $f$ is negative in the interior of
 $M$. From $\D f=-\lambda_1 f\geq 0$ on $M$ and $f|_{\pa M}=0$, we
 know from the the maximum principle and the Hopf lemma \cite{GT} that
 \be
 \fr{\pa f}{\pa \nu}>0 \ \ {\rm on}\ \  \pa M.
 \en
 It
 then follows from $0=\D f|_{\pa M}$ and (\ref{th5eq5}) that the mean
 curvature of $\pa M$ vanishes.
 Thus one can use Reilly theorem as stated before to conclude that $M$ is isometric to an $n$-dimensional  hemisphere of curvature $\kappa$.
 This completes the proof of Theorem \ref{th5}. \qed\\
\vs
 {\it Proof of Theorem \ref{th6}. } The discussions are similar to those in the proof of Theorem \ref{th5}. For the sake of completeness, we include it.
 Let $g$ be the eigenfunction of the problem (\ref{th6eq1}) corresponding to the first eigenvalue $\Lambda_1$:
\begin{equation}\label{th6eq3}
\left \{ \aligned &\Delta^2 g= -\Lambda_1 \D g   \ \ {\rm in \ \ } M, \\
&g= \fr{\pa^2 g}{\pa \nu^2}=0 \ \ {\rm on \ \ } \pa M.
\endaligned \right.
\end{equation}
Multiplying the first equation of (\ref{th6eq3}) by $g$ and
integrating on $M$, we have from  divergence theorem that \be
\label{th6eq4} \Lambda_1\int_M |\na g|^2= \int_M (\D g)^2 -\int_{\pa
M} s\D g, \en where $s=\left.\fr{\pa g}{\pa \nu}\right|_{\pa M}.$
Also, we have \be \label{th6eq5} \D g|_{\pa M} = (n-1) H s. \en
Hence \be \Lambda_1=\fr{\int_M (\D g)^2- (n-1)\int_{\pa M}
Hs^2}{\int_M|\na g|^2}, \en which, combining with Reilly's formula
gives \be \label{th6eq6} \Lambda_1&=&\fr{\int_M (|\na^2 g|^2+{\rm
Ric}(\na g, \na g)}{\int_M|\na g|^2}\\ \no &\geq & \fr{\int_M
\left(\fr{(\D g)^2}n +(n-1)\kappa |\na g|^2\right)}{\int_M|\na
g|^2}. \en We also have \be \int_M (\D g)^2 \geq \lambda_1\int_M
|\na g|^2,\en with equality holding if and only if $g$ is an
eigenfunction corresponding to $\lambda_1$. Therefore, we have
\be\label{th6eq7} \Lambda_1\geq \fr{\lambda_1}n +(n-1)\kappa.\en If
the equality occurs in (\ref{th6eq7}), then \be \D g =-\lambda_1 g \
\ {\rm on} \ \ M, \en and so \be \D^2 g= -\lambda_1 \D g.\en It then
follows that \be \Lambda_1=\lambda_1= \fr{\lambda_1}n
+(n-1)\kappa,\en which gives $0<\lambda_1=n\kappa$. As in the proof
of Theorem \ref{th5}, we can deduce that the mean curvature of $\pa
M$ vanishes and by Reilly's theorem that $M$ is isometric to an
$n$-dimensional Euclidean hemisphere of curvature $\kappa$.\qed\\

\section{Proof of Theorems \ref{th7}-\ref{th10} }
In this section, we shall prove Theorems \ref{th7}-\ref{th10}.

 {\it Proof of
Theorem \ref{th7}.} Let $w$ be an eigenfunction corresponding to the
first eigenvalue $q_1$  of the problem (\ref{eq1.6}):
\be\left\{\begin{array}{l}
 \Delta^2 w= 0  \ \ {\rm in \ \ } M, \\
w= \fr{\pa^2 w}{\pa \nu^2}  -q_1\fr{\pa w}{\pa \nu} =0 \ \ {\rm on \
\ } \pa M.
\end{array}\right.
\en Note that $w$ is not a constant since  $w|_{\pa M}=0$. Set
$\eta=\pa_{\nu} w|_{\pa M}$; then $\eta\neq 0$. Otherwise, we would
deduce from
$$w|_{\pa M}=\na w|_{\pa M}=\left.\fr{\pa^2 w }{\pa \nu^2}\right|_{\pa M}=0$$
that $\Delta w|_{\pa M}=0$ and so $\Delta w=0$ on $M$ by the maximum
principle, which in turn implies that $w=0$. This is a
contradiction.

From $w|_{\pa M}=0$, one gets from  divergence theorem that
\be\label{th7eq2} \int_M \langle \na w, \na(\Delta w)\rangle =
-\int_M w\Delta^2 w=0. \en Thus \be\label{th7eq3} \int_{\pa M} \D
w\fr{\pa w}{\pa \nu}& =&\int_M\langle\na (\D w), \na w\rangle
+\int_M (\D w)^2\\ \no &=&\int_M (\D w)^2. \en Since
\be\label{th7eq4}
\D w|_{\pa M}&=&\fr{\pa^2 w}{\pa \nu^2}+(n-1)H\fr{\pa w}{\pa\nu}\\
\no &=& q_1\fr{\pa w}{\pa\nu}+(n-1)H\fr{\pa w}{\pa\nu}, \en we
conclude that \be q_1=\fr{\int_M(\D w)^2-(n-1)\int_{\pa M}
H\eta^2}{\int_{\pa M} \eta^2}, \en which, combining with  Reilly's
formula, gives
 \be\label{th7eq5} q_1=\fr{\int_M(|\na^2 w|^2+{\rm
Ric}(\na w, \na w))}{\int_{\pa M} \eta^2} \geq \fr{\int_M |\na^2
w|^2}{\int_{\pa M} \eta^2}. \en The Schwarz inequality implies that
\be |\na^2 w|^2\geq \fr 1n(\D w)^2, \en with equality holding if and
only if \be\label{th7eq6} \na^2 w=\fr{\D w}n \langle, \rangle. \en
 Therefore, we have \be\label{th7eq7} q_1\geq \fr 1n\cdot
\fr{\int_M (\D w)^2}{\int_{\pa M} \eta^2}. \en On the other hand, we
have from the variational characterization of $p_1$ (cf.
(\ref{eq1.9})) that \be\label{th7eq8} \fr{\int_M (\D w)^2}{\int_{\pa
M} \eta^2}\geq p_1\en with equality holding if and only if $w$ is an
eigenfunction corresponding to $p_1$. Combining (\ref{th7eq7}) and
(\ref{th7eq8}), we get \be q_1\geq \fr{p_1}n.\en This proves the
first part of Theorem \ref{th7}.

Assume now that $q_1 = p_1/n$. In this case, we know from the above
proof that (\ref{th7eq6}) holds on $M$ and $w$ is an eigenfunction
corresponding to $q_1$. Hence
\be\label{th7eq9} \D w = p_1 \eta \ \
{\rm on } \ \ \pa M.\en Take an orthornormal frame $\{e_1,\cdots,
e_{n-1}, e_{n}\}$ on $M$ such that when restricted to $\pa M$,
$e_{n}=\nu$.
 From
$0=\na^2 w(e_i, e_{n})$, $i=1, \cdots, n-1$, and $w|_{\pa M}=0$, we
conclude that $\eta=b_0=const.\neq 0.$  It then follows from
 the
harmonicity of $\D w$ and $\ref{th7eq9})$ that \be\label{th7eq10} \D
w = p_1 b_0=nq_1 b_0 \ \ {\rm on\ \ } M. \en Setting
\be\label{th7eq10} \rho=-\fr{w}{p_1b_0}, \en we have
\be\label{th7eq11} \na^2 \rho = -\fr 1n \langle, \rangle \ \ {\rm
on\ \ } M. \en Taking the covariant derivative of (\ref{th7eq11}),
we get
 $\na^3 \rho=0$ and from the Ricci identity,
\be \label{th7eq12} R(X, Y)\na \rho=0, \en for any tangent vectors
$X, Y$ on $M$, where $R$ is the curvature tensor of $M$. By the the
maximum principle $\rho $ attains its maximum at some point $x_0$ in
the interior of $M$. Let $r$ be the distance function to $x_0$; then
from (\ref{th7eq11}) it follows that \be\label{th7eq13} \na \rho =
-\fr 1n r \fr{\pa}{\pa r}. \en Using (\ref{th7eq12}),
(\ref{th7eq13}), Cartan's theorem (cf. \cite{doC}) and $\rho|_{\pa
M}=0$, we conclude that $M$ is a ball in $\mathbb{R}^n$ whose center
is $x_0$, and \be \no \rho(x)=\fr 1{2n}(r_0^2 - |x-x_0|)\en in $M$,
here $r_0$ is the radius of the ball. Thus,
$$w(x)= \fr{p_1b_0}{2n } (|x-x_0|-r_0^2).$$
This completes the proof of Theorem \ref{th7}\qed\\

{\it Proof of Theorem \ref{th8} } Let $\phi$ be an eigenfunction
corresponding to the first eigenvalue $p_1$ of the problem
(\ref{eq1.5}), that is \be\left\{\begin{array}{l}
 \Delta^2 \phi= 0  \ \ {\rm in \ \ } M, \\
\phi= \D \phi -p_1 \fr{\pa \phi}{\pa \nu}=0  \ \ {\rm on \ \ } \pa
M.
\end{array}\right.
\en Set $\eta=\left.\fr{\pa \phi}{\pa \nu}\right|_{\pa M}$; then \be
\label{th8eq2}
 q_1=\fr{\int_{M}(\D \phi)^2}{\int_{\pa M} \eta^2}.
\en Substituting $\phi$ into  Reilly's formula, we have \be
\label{th8eq3}\int_{M} \left\{ (\D \phi)^2- |\na^2 \phi|^2\right\}&=
&\int_{M} {\rm Ric}(\na \phi, \na \phi)+\int_{\pa M} (n-1)H \eta^2\\
\no &\geq & -(n-1)\kappa \int_{M} |\na \phi|^2 + (n-1)c\int_{\pa M}
\eta^2\\ \no &\geq& -\fr{(n-1)\kappa}{\lambda_1} \int_{ M} (\D
\phi)^2 + (n-1)c\int_{\pa M}\eta^2, \en where, in the last step
above, we have used the fact that $\kappa\geq 0$ and
\be\label{th8eq4} \int_{M} |\na \phi|^2\leq \fr 1{\lambda_1}
\int_{M} (\D \phi)^2.\en The Schwarz inequality implies that
\be\label{th8eq5} |\na^2 \phi|^2\geq \fr 1n(\D \phi)^2 \en with
equality holding if and only if \be \label{th8eq6}\na^2 \phi=\fr{\D
\phi}n \langle, \rangle.\en Combining (\ref{th8eq2}), (\ref{th8eq3})
and (\ref{th8eq5}), we have \be\label{th8eq7} p_1\geq\fr{
nc\lambda_1}{n\kappa+\lambda_1}.\en If the equality sign holds in
(\ref{th8eq7}), then the inequalities (\ref{th8eq3}) and
(\ref{th8eq5}) must take equality sign and, in particular,
(\ref{th8eq6}) holds on $M$. By using the same arguments as in the
proof of Theorem \ref{th7}, we conclude that $\D\phi$ is a nonzero
constant and $M$ is isometric to a ball in $\mathbb{R}^n $ which, in
turn, implies that $\kappa =0$.
\qed\\

{\it Proof of Theorem \ref{th9}.}  Let $u$ be an eigenfunction of
the problem (\ref{eq1.8}) with $M=\mathbb{B}$ corresponding to an
eigenvalue $\va$. From the bi-harmonicity of $u$ we can two uniquely
determined harmonic functions $g, h$ on $\mathbb{B}$, such that
\cite{AC} \be\label{th9eq1} u(x)= g(x)+|x|^2 h(x). \en
 Let us denote by
$\mathcal{E}$ the Euler operator defined by\be\label{th9eq2} \ma
f(x):=\sum_{i=1}^n x_i\fr{\pa f}{\pa x_i}(x). \en
 We will use the
notation $\ma^lf=\ma (\ma^{l-1}f), l=1, 2,\cdots. $ Observe that
$\ma$ is related to the exterior normal derivative by the relation
\be\label{th9eq3} \left.\fr{\pa f}{\pa\nu}\right|_{\pa\mathbb{
B}}=\ma f|_{\pa\mathbb{B}},\en and if $f$ is harmonic then so is
$\ma f$. Setting $-2w=g+h$, we have $u(x)=-2w(x)+(|x|^2-1)h(x)$.
From
$$0=\pa_{\nu} u|_{\pa\mathbb{B}}=(-2\pa_{\nu} w+ 2h)|_{\pa \mathbb{B}}=(-2\ma w+ 2h)|_{\pa B}$$
and the harmonicity of  $-2\ma w+ 2h$, we know that $-2\ma w+ 2h=0$
on $\mathbb{B}$. Hence \be\label{th9eq4} u(x)=-2w(x)+(|x|^2-1)\ma
w(x).\en We have
 \be\label{th9eq5} \D u= 2n\ma w +4\ma^2w, \en and
so \be\label{th9eq6} \pa_{\nu}(\D u)|_{\pa\mathbb{B}}=(2n \ma^2w
+4\ma^3w)|_{\pa\mathbb{B}}. \en If $z\in\pa\mathbb{B}$, we have from
$\pa_{\nu} u(z)=0$ that \be\label{th9eq7} \D u(z)=\ov{\D} u(z)+\na^2
u(\nu, \nu)(z)=\ov{\D} u(z)+\ma^2 u(z). \en For any
$x\in\mathbb{B}$, a simple calculation gives \be\label{th9eq8} \ma
u(x)=(|x|^2-1)(2\ma w +\ma^2 w)(x) \en and \be \label{th9eq9} \ma^2
u(x)= 2|x|^2 (2\ma w +\ma^2 w)(x)+(|x|^2-1)(2\ma^2 w +\ma^3 w)(x).
\en Thus \be \ma^2u(y)= 2 (2\ma w +\ma^2 w)(y), \ \ \forall y\in \pa
\mathbb{B}, \en which, combining with (\ref{th9eq5}) and
(\ref{th9eq7}), gives \be\label{th9eq10} \ov{\D} u(y)=((2n-4)\ma
w+2\ma^2 w)(y), \ \ \forall y\in \pa\mathbb{B}. \en Substituting
(\ref{th9eq4}), (\ref{th9eq6}) and (\ref{th9eq10}) into the
equation: \be\no \pa_{\nu}(\D u) + \beta \overline{\D} u + \varsigma
u =0 \ \ {\rm on \ \ } \pa\mathbb{B}, \en we infer
\be\label{th9eq11}
 0= (2n\ma^2w +4\ma^3w
+\beta((2n-4)\ma w+2\ma^2 w)-2\varsigma w)|_{\pa\mathbb{B}}. \en It
follows that \be \label{th9eq12} n\ma^2w +2\ma^3w +\beta((n-2)\ma
w+\ma^2 w)=\varsigma w \ \ {\rm on} \ \ \mathbb{B}. \en Since $w$ is
harmonic on $\mathbb{B}$, there exist $p_m\in {\mathcal D}_m$ such
that \be\label{th9eq13} w(x)=\sum_{m=0}^{\infty} p_m(x) \en for all
$x\in B$, the series converging absolutely and uniformly on compact
subsets of $\mathbb{B}$(cf. \cite{ABR}). Substituting
(\ref{th9eq13}) into (\ref{th9eq12}) and using $\ma p_m = mp_m$, we
infer \be\label{th9eq14} 2m^3+ nm^2+\beta(m^2+(n-2)m) = \va, m=0,
1,\cdots \en which shows that all but one of the $p_{m}'s$ are
zeros. Consequently, the eigenvalues are \be \label{th9eq15}
\va_k=k^2(n+2k)+\beta k(k+n-2), k=0, 1, \cdots,\en the multiplicity
of $\va_k$ is the dimension of ${\mathcal D}_k$ and the eigenspace
corresponding to $\va_k$ is \be \mathbf{E}_k=\{
-2w_k+k(|x|^2-1)w_k|\ \  w_k\in {\mathcal D}_k\}. \en This completes
the proof of Theorem \ref{th9}. \qed\\

 {\it Proof of Theorem \ref{th10}.}
From (\ref{eqxi}), we only need to show that the first non-zero
eigenvalue $\xi_1$ of the problem (\ref{eq1.7}) satisfies
\be\label{th10eq2} \xi_1 >\fr{nc\lambda_1 \mu_1}{(n-1)(\mu_1+n\kappa
)}. \en
 Let $f$ be an eigenfunction corresponding $\xi_1$:
 \be\label{th10eq3}
\left\{\begin{array}{l} \Delta^2 f= 0  \ \ {\rm in \ \ } M, \\
\fr{\pa }{\pa \nu} f = \fr{\pa }{\pa \nu} \D f+\xi_1 f =0 \ \ {\rm
on \ \ } \pa M.
\end{array}
\right. \en Set $z=f|_{\pa M}$; then $z\neq 0$ and
\be\label{th10eq4} \xi_1=\fr{\int_{M}(\D f)^2}{\int_{\pa M} z^2}.
\en Substituting $f$ into  Reilly's formula, we have \be
\label{th10eq5}\int_{M} \left\{ (\D f)^2- |\na^2 f|^2\right\}&=&
\int_M {\rm Ric}(\na f, \na f)+ \int_{\pa M} II(\overline{\na}z,
\ov{\na} z)
\\ \no &\geq & -(n-1)\kappa \int_M |\na f|^2 + c\int_{\pa M} |\ov{\na}z|^2.
\en Since $\pa_{\nu} f|_{\pa M}=0$, we have \be\label{th10eq6}
\int_M (\D f)^2\geq \mu_1 \int_M |\na f|^2. \en It follows from
(\ref{th10eq3}) that $\int_{\pa M} z=0$ and so we have from the
Poincar\'e inequality that \be\label{th10eq7} \int_{\pa
M}|\ov{\na}z|^2 \geq \lambda_1 \int_{\pa M} z^2.\en The Schwarz
inequality implies that \be\label{th10eq8} |\na^2 f|^2\geq \fr 1n(\D
f)^2 \en with equality holding if and only if $\na^2 f=\fr{\D f}n
\langle, \rangle$.

 Combining (\ref{th10eq4})-(\ref{th10eq8}),  we get
\be\label{th10eq9}\xi\geq \fr{nc\lambda_1 \mu_1}{(n-1)(\mu_1+n\kappa
)}. \en Let us show by contradiction that the equality  in
(\ref{th10eq9})  can't occur. In fact, if (\ref{th10eq9}) take
equality sign,  then we must have \be\label{th10eq10} \na^2 f=\fr{\D
f}n \langle, \rangle. \ \ {\rm on \ \ } M.\en Thus for a tangent
vector field $X$ of $\pa M$,  we have from (\ref{th10eq10}) and
$\pa_{\nu}f|_{\pa M}=0$ that \be 0=\na^2f(\nu, X)=X\nu
f-(\na_X\nu)f=-\langle \na_X \nu, \overline{\na} z\rangle. \en In
particular, we have \be\no II(\overline{\na }z, \overline{\na }z)=0.
\en
This is impossible since $II\geq c I$ and $ z$ is not  constant. This finishes the proof Theorem \ref{th10}.\qed \\

An interesting question related to the above proof is to find a
sharp lower bound for $\xi_1$. We propose a \vs {\bf Conjecture.}
{\it  Let $M$ be an $n$-dimensional compact Riemannian manifold with
boundary and nonnegative Ricci curvature. Assume that the principal
curvatures of $\pa M$ are bounded below by a positive constant $c$
and denote by $\lambda_1$ the first nonzero eigenvalue of the
Laplacian of $\pa M$. Then the first nonzero eigenvalue of the
problem (\ref{eq1.7}) satisfies \be\label{coneq1} \xi_1\geq
\fr{(n+2)c\lambda_1}{n-1}\en with equality holding if and only if
$M$ is isometric to an $n$-dimensional Euclidean ball of radius
$1/c$.}

\begin{flushleft}

Qiaoling Wang\\
Departamento de Matem\'atica\\
Universidade de
Bras\'{\i}lia\\
Bras\'{\i}lia-DF 70910-900, Brazil\\
e-mail: wang@mat.unb.br
\end{flushleft}
\vskip0.5cm
\begin{flushleft}
Changyu Xia\\
Departamento de Matem\'atica\\
Universidade de
Bras\'{\i}lia\\
Bras\'{\i}lia-DF 70910-900 , Brazil\\
e-mail: xia@mat.unb.br
\end{flushleft}

\end{document}